\documentclass[12pt]{article}

\usepackage{amssymb}
\usepackage{colortbl}
\usepackage{amsmath}

\usepackage{moreverb}
\usepackage{booktabs}

\usepackage{amsthm}

\newtheorem{thm}{Theorem}
\newtheorem{cor}{Corollary}
\newtheorem{defi}{Definition}

\newcommand{\norm}{|\hspace{-1pt}|\hspace{-1pt}| }

\begin{document}


\title{Multiple summation inequalities and their application to stability analysis of discrete-time delay systems}

\author{\'E. Gyurkovics\footnote{Corresponding  author} \footnote{Mathematical Institute, Budapest University of Technology and Economics; gye@math.bme.hu, kk@math.bme.hu, ilona.nagy@gmail.com }, \hspace{0.2cm}  K. Kiss\footnotemark[\value{footnote}], \hspace{0.2cm}
  I. Nagy\footnotemark[\value{footnote}], \hspace{0.2cm}  T. Tak\'acs\footnote{Corvinus University of Budapest; Tibor.Takacs@kih.gov.hu}}


\date{}

\maketitle

\begin{abstract}
This paper is devoted to stability analysis of discrete-time delay systems based on a set of Lyapunov-Krasovskii functionals. New multiple summation inequalities are derived that involve the famous discrete Jensen's and Wirtinger's inequalities, as well as the recently presented inequalities for single and double summation in \cite{Nam15b}. The present paper aims at showing that the proposed set of sufficient stability conditions can be arranged into a bidirectional hierarchy of LMIs establishing a rigorous theoretical basis for comparison of conservatism of the investigated methods. Numerical examples illustrate the efficiency of the method.
\end{abstract}



\section{Introduction}
Time delays are frequently encountered in various real life phenomena, e.g. in
physical, industrial, and engineering systems.
Since time delays
 may result in instability and poor performance of
the systems,  the stability of systems with time delays
has received much attention during the past few decades:
a comprehensive review can be found e.g. in \cite{Briat14}, \cite{frid2014}, \cite{WHSh}.
(See also the references therein.)

This paper investigates the stability issue of linear discrete-time delay systems described by
\begin{eqnarray}
x(t+1)&=&Ax(t)+A_d x(t- \tau ), \hspace{0.5cm}  t=0,1,...  \label{uj9} \\
x(s)&=& x_0 (s), \hspace{3.0cm} s=-\tau, -\tau+1, ... ,0, \nonumber
\end{eqnarray}
where $x(t) \in R^{n_x}$ is the state, $A$ and $A_d$ are given constant matrices of appropriate size, the time delay $\tau $ is a known positive integer and $x_0(.)$ is the initial function.
It is well-known that a necessary and sufficient condition can be derived for the stability of (\ref{uj9}) by the so called lifting technique in the form of the discrete Lyapunov inequality. This approach, however, suffers from the curse of dimensionality, since the number of decision variables in the LMI to be solved is $(\tau +1)n_x((\tau+1) n_x +1)/2,$ which may be too large, if $\tau$ is large. Therefore, much effort has been devoted to obtain sufficient conditions that are computationally less intensive. The Lyapunov-Krasovskii functional (LKF) approach offers a fruitful alternative: during the past decades numerous LKFs have been proposed, and simultaneously, a large number of different techniques has been developed to get better estimations for the forward difference of the functional. Jensen's inequality (a genaral form and historical comments can be found in \cite{marek}, its role in delayed systems is apparent from references \cite{Briat14}, \cite{frid2014}, \cite{gyt15}, \cite{hien15}, \cite{kwon14}-\cite{WHSh}, \cite{ZhH}) and  Wirtinger's inequality (see \cite{Hardy} for original version and \cite{gye15}, \cite{hien15}, \cite{kwon14}-\cite{LFrAut12}, \cite{Nam15a}-\cite{seur15}, \cite{ZhH}, \cite{zha15} for generalizations) play an outstanding role in this respect. Recently, efficient inequalities for single and double summation have been published in \cite{Nam15b} that involves the previous two inequalities. In this way, the derived sufficient conditions yield better and better results in respect of the tolerated delay bound. (See e.g. \cite{hien15}, \cite{kwon14}-\cite{seur15}, \cite{ze15}-\cite{zha15} and the references therein.) In the vast majority of cases, however, the superiority of methods is demonstrated by comparison of the results obtained for some benchmark examples.

The excellent idea of a hierarchy of sufficient LMI  stability conditions has been introduced in \cite{seur14, seur14b} that gives - among the others - a rigorous basis for comparison.

\emph{The aim of the present work is twofold. Firstly, new inequalities will be derived to the case of arbitrary number of summation both for functions and differences. Secondly, a set of LKFs will be proposed, and a hierarchy table of the obtained sufficient stability conditions will be demonstrated.
We shall show by some benchmark examples that the proposed methods give better upper bounds for the tolerable time delay than the best ones that we could find in the previously published literature.}

The notations applied in the paper are very standard, therefore we mention only a few of them. $\mathbf{P} _K$ denotes the space of polynomials having degree not higher than $K$. Symbol $A \otimes B$ denotes the Kronecker-product of matrices $A,B,$ while $\mathbf{S}_{n}^{+}$ is the set of positive definite symmetric matrices of size $n\times n.$

\section{Multiple summation inequalities}
\subsection{Discrete orthogonal polynomials}
Suppose that $m$ and $N$ are given positive integers, and consider the support points $s_i=i,$ $(i=0,1,\ldots,N-1)$. 
For functions $f,g: \mathbf{Z} \rightarrow \mathbf{R},$ define a scalar product by
\begin{eqnarray}
\ll f,g \gg _{m} &=& \sum _{i_1 =0}^{N-1} \sum _{i_2 =0}^{i_1} ... \sum _{i_m =0}^{i_{m-1}} f(i_m)g(i_m),      \label{eg3}
\end{eqnarray}
and denote the corresponding norm by $\norm f \norm _m.$
Since
\begin{eqnarray*}
\sum _{j_1 =0}^{K} \sum _{j_2 =0}^{j_1} ... \sum _{j_{m+1} =0}^{j_{m}} 1= \left(
                                                                            \begin{array}{c}
                                                                              K+m+1 \\
                                                                              m+1 \\
                                                                            \end{array}
                                                                          \right),
\end{eqnarray*}
the scalar product in (\ref{eg3}) can equivalently be written as
\begin{eqnarray}
 \ll f,g \gg _{m}
 & =&
  \sum _{i_m =0}^{N-1} f(i_m)g(i_m) \sum _{i_1 = i_m }^{N-1} ... \sum _{i_{m-1} = i_m}^{i_{m-2}}1
  \nonumber\\
   & = &\frac{1}{(m-1)!} \sum _{i =0}^{N-1} r_{N,m-1}(i)f(i)g(i). \label{eg3a}
\end{eqnarray}
where $r_{N,m}(i)=m! \left(
                                                                             \begin{array}{c}
                                                                               N-1+m-i \\
                                                                               m \\
                                                                             \end{array}
                                                                           \right)$,
$i=0,1,...,N-1$.
For convenience, we introduce another scalar product by
\begin{equation}
<f,g>_m = \sum _{i =0}^{N-1}  r_{N,m-1}(i)f(i)g(i), \hspace{0.5cm}    \label{eg4}
\end{equation}
and denote the corresponding norm by $\left\| f \right\|_{m} .$ Obviously,
\begin{equation}
  \| f \| ^{2}_{m} = (m-1)! \ \norm f \norm _{m}^{2}. \label{eg5}
\end{equation}
Let $p_{ml}(.)$ denote the discrete
orthogonal polynomials on 
 $\mathbf{Z} \cap [0,N-1]$ with respect to the the scalar product (\ref{eg3})(or (\ref{eg4})), and with the exact
degree $l=0,...,N-1.$
It is well-known that these polynomials can be generated by applying either the Gram - Schmidt orthogonalization process or a three term recurrence relation (see e.g. \cite{gau}). If $m=1$, the corresponding orthogonal polynomials are called discrete Chebyshev polynomials, and they are given by the three term recurrence relation together with their norms e.g. in \cite{gau}. However, we are not aware of similar published formulas for $m>1.$ Therefore the polynomials $p_{mj}$ have been generated for $m>1$ using Wolfram Mathematica and the results can be found in \cite{kisnagy}.
%
\subsection{Summation inequalities for functions}
Let $R \in \mathbf{S}^{+}_{n}$ be given. 
For any $f: \mathbf{Z} 
\rightarrow \mathbf{R}^{n},$
consider the functional
\begin{eqnarray*}
J_m (f) = \sum _{i_1 =0}^{N-1} \sum _{i_2 =0}^{i_1} ... \sum _{i_m =0}^{i_{m-1}} f^{T}(i_m)R f(i_m).
\end{eqnarray*}

\emph{Our aim in this subsection is to derive lower bounds for functionals of this type.}



Let $\nu _1, \ \nu _m$ be  given nonnegative integers satisfying condition $\nu_m < \nu_1<N.$ Set $\widetilde{w}_{mj}=\ll f,p_{mj} \gg _{m}$ and $w_{mj}= <f,p_{mj}>_{m}$, where $j=0,\ldots, \nu _m,$ and the scalar product is taken componentwise. Clearly,
$w_{mj} = (m-1)! \ \widetilde{w}_{mj}.$ In what follows, $w_{1j}$ will play a key role, therefore we shall denote it specially
\begin{equation}\label{gy1}
 \phi _j = w_{1j}, \; \; j=0,1,\ldots, \nu _1.
\end{equation}

Suppose that $\nu _m +m-1 \leq \nu _1$ and set $q_{N,m+j-1}(i)=r_{N,m-1} (i)p_{mj}(i).$
Then $q_{N,m+j-1} \in \mathbf{P} _{\nu _1}$  with degree $m+j-1 \leq \nu _m +m-1 < \nu _1$.
Since the polynomials $p_{1l}$, $l=0,...,\nu _1$  form a basis of $\mathbf{P} _{\nu _1},$
we have
\begin{eqnarray}
&&q_{N,m+j-1}(i)=\sum _{l=0}^{m+j-1}\xi _{j,l}^{m}p_{1l}(i).    \label{uj1}
\end{eqnarray}
Introduce now the row vector of length $\nu _1+1$
\begin{eqnarray*}
&&\underline{\xi} _{j}^{m} =  \left[ \xi _{j,0}^{m} \; \xi _{j,1}^{m} \; \ldots \; \xi _{j,\nu _1}^{m}  \right]
\end{eqnarray*}
the first $j+m$ entries of which equal to the coefficients in (\ref{uj1}), while 
$\xi _{j,l}^{m}=0$ if $j+m-1 < l \leq \nu _1$.

 Then we have
\begin{eqnarray}
&&w_{mj}=\sum _{l=0}^{\nu _1} \xi _{j,l}^{m}\sum _{i=0}^{N-1}p_{1l}(i)f(i)=\sum _{l=0}^{\nu _1} \xi _{j,l}^{m} \phi _l.\label{gy3}
\end{eqnarray}
Introduce the notations:
\begin{eqnarray}
&& \Xi_m =     \left(
                    \begin{array}{c}
                      \underline{\xi} _{0}^{m} \\
                      \vdots \\
                      \underline{\xi} _{\nu _m}^{m} \\
                    \end{array}
                  \right)
    \in  \mathbf{R}^{(\nu _m +1) \times (\nu _1 +1)}, \hspace{1cm} \Phi= \mbox{col} \left\{ \phi _0, \; \ldots ,\; \phi _{\nu_1} \right\},
    \label{gy4} \\
&&  \mathcal{R}_m = \mbox{diag} \left\{ \chi _{m,0}R, \; \ldots , \;
\chi _{m,\nu _m}R  \right\},  \hspace{0.8cm} \chi _{m,j}=\frac{1}{\|p_{mj}\|_{m}^{2}}.\label{gy5}  
\end{eqnarray}

\begin{thm} \label{Th:1}
 Let  $m, \nu _1, \nu _m, N$  be given integers satisfying conditions $m\geq 1$ and $\nu_m < \nu_1<N.$  Let matrix $R \in \mathbf{S}^{+}_{n}$ and function $f: \mathbf{Z} 
\rightarrow \mathbf{R}^{n}$ be given. Then the following inequality holds:
\begin{eqnarray}
J_m (f)\geq \frac{1}{(m-1)!} \ \Phi ^T \left( \Xi _m \otimes I \right)^T \mathcal{R}_m \left( \Xi _m \otimes I \right) \Phi.  \label{eg6}
\end{eqnarray}
\end{thm}
\textbf{Proof.}
Let $\mu _j$ denote arbitrary constants and set
\begin{eqnarray*}
z(i)=f(i)-\sum _{j =0}^{\nu _m} \mu _j \widetilde{w}_{mj} p_{mj} (i).
\end{eqnarray*}
Then we obtain
\begin{equation}
0 \leq J_m (z)= J_m (f)- \sum _{j=0}^{\nu_m} (2 \mu _j - \mu_{j}^{2} \norm p_{mk} \norm _{m}^{2} )
 \widetilde{w}_{mj}^{T} R \widetilde{w}_{mj}.        \label{eg6a}
\end{equation}
The term $2 \mu _j - \mu_{j}^{2} \norm p_{mk} \norm _{m}^{2}$ takes its maximum at
$\mu _j^{\ast} = \frac{1}{\norm p_{mk} \norm _{m}^{2}}$. By rearranging (\ref{eg6a}) and by substituting $\mu _j^{\ast}$  one obtains
 that
 \begin{eqnarray}
J_m (f) &\geq& \sum _{j =0}^{\nu _m} \frac{1}{\norm p_{mk} \norm _{m}^{2}} \widetilde{w}^{T}_{mj} R \widetilde{w}_{mj}. \label{gy2}
\end{eqnarray}
Now, one has only to take into account (\ref{eg5}), (\ref{gy3})-(\ref{gy5}) and the relation between $\widetilde{w}_{mj}$ and ${w}_{mj}$ to obtain (\ref{eg6}).
  $\Box$

From Theorem \ref{Th:1}, one can easily derive several known summation inequalities such as the discrete Jensen inequality (see e.g. ( \cite{gyt15}, \cite{hien15}, \cite{marek},  \cite{lee15}, \cite{Nam15a}, \cite{Nam15b}, \cite{seur15}, \cite{ZhH}), the discrete Wirtinger inequality (\cite{hien15}, \cite{lee15}, \cite{Nam15a}, \cite{Nam15b}), and the inequalities of Nam, Trinh and Pathirana \cite{Nam15b}.
As it is usual in the literature of stability analysis (see e.g. \cite{hien15}), the lower estimations of type (\ref{eg6}) are given by
the single, double, triple, etc. summation of the state variable. We follow this line too, when formulating the corollaries of
Theorem 1.

\begin{cor} \label{cor:1}
If $m=1$,and $f(i)=x(i),$  (\ref{eg6}) implies
\begin{itemize}
\item
 for $\nu _1=0$ the Jensen inequality
\begin{equation}
J_1 (x) \geq  \frac{1}{N}  \Omega _{10}^T
R \Omega _{10} \hspace{1cm} \mbox{with} \hspace{0.5cm} \Omega _{10} = \sum _{i_1 =0}^{N-1} x(i_1) ;    \label{eg7}
\end{equation}
\item
for $\nu _1 =1$, $N>1$, the Wirtinger inequality
\begin{eqnarray}
J_1 (x) &\geq &\frac{1}{N}\left\{ \Omega_{10}^{T} R \Omega _{10} +      
3 \frac{N+1}{N-1} \Omega_{11}^{T} R \Omega _{11} \right\}    \label{eg8} \\ 
&& \hspace{-5.65cm}\mbox{with} \nonumber  \\
\Omega _{11} &=& \sum _{i_1 =0}^{N-1} x(i_1) - \frac{2}{N+1} \sum _{i_1 =0}^{N-1} \sum _{i_2 =0}^{i_1} x(i_2)
; \label{gy10}
\end{eqnarray}
\item
for $\nu _1 =2$, $N>1$,  the inequality
\begin{eqnarray}
&& \hspace*{-0.7cm} J_1 (x) \geq \frac{1}{N}\left\{ \Omega_{10}^{T} R \Omega _{10} +      
3 \frac{N+1}{N-1} \Omega_{11}^{T} R \Omega _{11} + 5 \frac{(N+1)(N+2)}{(N-1)(N-2)} \Omega_{12}^{T} R \Omega _{12}
\right\} \label{eg11} \\
&&\hspace{-1.6cm}\mbox{with} \nonumber  \\
&& \hspace*{-0.7cm} \Omega _{12} = \sum _{i_1 =0}^{N-1} x(i_1)-\frac{6}{N+1} \sum _{i_1 =0}^{N-1} \sum _{i_2 =0}^{i_1} x(i_2)
+\frac{12}{(N+1)(N+2)} \sum _{i_1 =0}^{N-1} \sum _{i_2 =0}^{i_1} 
\sum _{i_3 =0}^{i_2} x(i_3). \label{gy11}
\end{eqnarray}
\end{itemize}
\end{cor}

\textbf{Proof.}
Since the first three discrete Chebishev polynomials and their norms are known to be
\begin{eqnarray}
&&\hspace*{-0.5cm}p_{10} (x) \equiv 1, \hspace{6.5cm} \| p_{10} \| ^{2}_{1} =N, \label{cheb0}\\
&&\hspace*{-0.5cm}p_{11} (x)= 2x+1-N, \hspace{4.7cm} \| p_{11} \| ^{2}_{1} = \frac{N(N^{2}-1)}{3},   \label{cheb1}\\
&&\hspace*{-0.5cm}p_{12} (x)= 6x^{2}-6(N-1)x+(N-1)(N-2),   \hspace{0.4cm}
\| p_{12} \| ^{2}_{1} = \frac{(N^{2}-4)(N^{2}-1)N}{5},    \label{cheb2}
\end{eqnarray}
(see e.g. \cite{gau}), the proof consists of some straightforward but lengthy computations the details of which are omitted. $\Box$

\textbf{Remark 1.} Estimation (\ref{eg11})-(\ref{gy11}) is less conservative than that in Lemma 1 of \cite{hien15},
and it is identical with that of \cite{Nam15b}, equation (27).

\begin{cor} \label{cor:2}
If $m=2$,and $f(i)=x(i),$  (\ref{eg6}) implies
\begin{itemize}
\item
 for $\nu _2=0$,  inequality
\begin{equation}
J_2 (x) \geq  \frac{2}{N(N+1)}  \Omega _{20}^T
R \Omega _{20} \hspace{1cm} \mbox{with} \hspace{0.5cm} \Omega _{20} = \sum _{i_1 =0}^{N-1} \sum _{i_2 =0}^{i_1} x(i_2) ;  \label{eg20}
\end{equation}
\item
for $\nu _2 =1$, $N>1$, inequality
\begin{eqnarray}
J_2 (x) &\geq &\frac{2}{N(N+1)}\left\{ \Omega_{20}^{T} R  \Omega _{0} +      
8 \frac{N+2}{N-1} \Omega_{21}^{T} R \Omega _{21} \right\}     \label{eg22} \\ 
&& \hspace{-4.8cm}\mbox{with} \nonumber  \\
\Omega _{21} &=& \sum _{i_1 =0}^{N-1} \sum _{i_2 =0}^{i_1} x(i_2) -
\frac{3}{N+2} \sum _{i_1 =0}^{N-1} \sum _{i_2 =0}^{i_1} \sum _{i_3 =0}^{i_2} x(i_3). \label{gy20}
\end{eqnarray}
\end{itemize}
\end{cor}
\textbf{Proof.}
Applying the first two orthogonal  polynomials
and their norms
\begin{eqnarray*}
&&p_{20} (x) \equiv 1, \hspace{2.5cm} \| p_{20} \| ^{2}_{2} =\frac{N(N+1)}{2},  \\
&&p_{21} (x)= x+\frac{1-N}{3}, \hspace{0.8cm}
\| p_{21} \| ^{2}_{2} = \frac{(N-1)N(N+1)(N+2)}{36},
\end{eqnarray*}
the proof consists of some straightforward but lengthy computations the details of which are omitted. $\Box$

\textbf{Remark 2.} Estimation (\ref{eg22}) is equivalent to that of Lemma 2 in \cite{hien15}.

\subsection{Summation inequalities for differences}
Based on Theorem \ref{Th:1}, we want to derive lower estimations for the functional $J_m$ applied to the forward difference of a function. In doing so, let  $\rho$ be given by function 
$f: \mathbf{Z} 
\rightarrow \mathbf{R}^n$ as
\begin{equation}\label{gy30}
  \rho : \mathbf{Z} 
  \rightarrow \mathbf{R}^n, \; \rho(i) =f(i+1)-f(i),\; i=0,1,...,N-1.
\end{equation}

\begin{thm}\label{Th:2}
 Let  $m, \nu _1, \nu _m, N$ be given integers satisfying condition $\nu _m +m-1 \leq \nu _1<N$. Let  $R \in \mathbf{S}^{+}_{n}$ be given and let function $\rho$ defined by (\ref{gy30}). Then
the following inequality holds:
\begin{eqnarray} \label{gy60}
J_m (\rho)\geq \frac{1}{(m-1)!} \widetilde{\Phi} ^T \left( \mathcal{Z} _m \otimes I \right)^T \mathcal{R}_m \left( \mathcal{Z} _m \otimes I \right)
\widetilde{\Phi} ,
\end{eqnarray}
where $\widetilde{\Phi}=\mbox{col}\left\{ f(N),f(0),\phi_0,\ldots ,\phi_{\nu_1-1} \right\}$, $\mathcal{R}_m$ is defined by (\ref{gy5}) and $\mathcal{Z} _m $ is given by (\ref{gy50}) below.
\end{thm}
\textbf{Proof.}
Introduce the notation $\theta _{mj} = <\rho , p_{mj}>_m,$ $j=0,...,m,$ where the scalar product is taken componentwise. 
From (\ref{eg6a}), it follows immediately that
\begin{equation}
 J_m(\rho) \geq \frac{1}{(m-1)!}\sum_{j=0}^{\nu_m}\frac{1}{\| p_{mj} \| _{m}^{2}} \theta _{mj}^{T} R \theta _{mj}.                 \label{eg25}
\end{equation}
 The vectors $\theta _{mj}$ have to be expressed by vectors $\phi _l = w_{1l}.$
\begin{eqnarray*}
\theta _{m,j} =  <\rho,p_{mj}>_m&=&\sum _{i=0}^{N-1}r_{N,m-1}(i)p_{mj}(i)f(i+1) - \sum _{i=0}^{N-1}r_{N,m-1}(i)p_{mj}(i)f(i) \\
 &=&
 r_{N,m-1}(N-1)p_{mj}(N-1)f(N)-r_{N,m-1}(-1)p_{mj}(-1)f(0)   \\
&&+\sum _{i=0}^{N-1}\left[r_{N,m-1}(i-1)p_{mj}(i-1) - r_{N,m-1}(i)p_{mj}(i)\right]f(i). 
\end{eqnarray*}
Apply now the notation
\begin{eqnarray}
&&\widetilde{q}_{N,m+j-2}(i)= r_{N,m-1}(i-1)p_{mj}(i-1)- r_{N,m-1}(i)p_{mj}(i).    \label{uj2}
\end{eqnarray}
The degree of polynomial (\ref{uj2}) is exactly $m+j-2$, and $m+j-2\leq m+ \nu _1 -2 \leq \nu _1 -1$, therefore there
exist coefficients $\zeta _{j,l}^{m}$ such that
\begin{eqnarray}
&& \widetilde{q}_{N,m+j-2}(i)=\sum _{l=0}^{m+j-2} \zeta _{j,l}^{m} p_{1l}(i).   \label{uj3}
\end{eqnarray}
Set  $c_{m,j,0}=-r_{N,m-1}(-1)p_{mj}(-1)$ and $c_{m,j,1}=r_{N,m-1}(N-1)p_{mj}(N-1).$
Then we have
\begin{eqnarray}\label{gy65}
&&\theta _{m,j} = c_{m,j,1}f(N)+c_{m,j,0}f(0)+\sum _{l=0}^{m+j-2} \zeta _{j,l}^{m} \phi_l .
\end{eqnarray}
Introduce now the row vector of length $\nu _1 +2$
\begin{eqnarray*}
&&\underline{\zeta} _{j}^{m} =  \left[ c_{m,j,1} \; c_{m,j,0} \; \zeta _{j,0}^{m}
\; \ldots \; \zeta _{j,\nu _1-1}^{m} \right] ,
\end{eqnarray*}
where entries $\zeta _{j,l}^{m}$ equal to the corresponding coefficients of (\ref{uj3}) 
for $0\leq l\leq m+j-2$, and $\zeta _{j,l}^{m}=0$ is taken for $m+j-2 < l \leq \nu _1-1$.

We note that $\zeta _{\nu_m,m+\nu_m-2}^{m}\neq 0,$ since the exact degree of $\widetilde{q}_{N,m+\nu_m-2}$ is $m+\nu_m-2,$ and
$\zeta _{\nu_m,m+\nu_m-2}^{m}$ is the coefficient of the basis element $p_{1,m+\nu_m-2}.$
Set
\begin{eqnarray} \label{gy50}
&&\mathcal{Z}_m = \left(
                    \begin{array}{c}
                      \underline{\zeta}_{0}^{m} \\
                      \vdots \\
                      \underline{\zeta}_{\nu_m}^{m} \\
                    \end{array}
                  \right) \in  \mathbf{R}^{(\nu _m +1) \times (\nu _1 +2)}.
\end{eqnarray}
Then (\ref{gy60}) immediately follows from (\ref{eg25}), (\ref{gy65}) and (\ref{gy50}). $\Box$

\begin{cor} \label{cor:3}
 If $m=1$, and $f(i)=x(i),$  (\ref{gy60}) implies
\begin{itemize}
\item
 for $\nu _1=0$ the discrete Jensen inequality for differences
\begin{equation}
J_1 (\rho) \geq  \frac{1}{N}  \widetilde{\Omega} _{10}^T
R \widetilde{\Omega} _{10} \hspace{1cm} \mbox{with} \hspace{0.5cm} \widetilde{\Omega} _{10} =  x(N)-x(0) ;    \label{eg7uj}
\end{equation}
\item
for $\nu _1 =1$, $N>1$, the discrete Wirtinger inequality for differences
\begin{eqnarray}
J_1 (\rho) &\geq &\frac{1}{N}\left\{ \widetilde{\Omega}_{10}^{T} R \widetilde{\Omega} _{10} +      
3 \frac{N+1}{N-1} \widetilde{\Omega}_{11}^{T} R \widetilde{\Omega} _{11} \right\} ,   \label{eg8uj} \\ 
&& \hspace{-5.65cm}\mbox{with} \nonumber  \\
\widetilde{\Omega} _{11} &=& x(N)+x(0)-\frac{2}{N+1} \sum _{i_1 =0}^{N} x(i_1) .\label{gy10uj}
\end{eqnarray}
\end{itemize}
\end{cor}
\textbf{Proof.} We calculate with the polynomials (\ref{cheb0})-(\ref{cheb1}) again: (\ref{eg7uj}) is evident, while (\ref{eg8uj})-(\ref{gy10uj}) follows by taking into account that
\begin{eqnarray*}
 \theta_{1,1} &=& (N-1)x(N)+(N+1)x(0)-2 \sum _{i_1 =0}^{N-1} x(i_1). \hspace{1cm}\Box
 \end{eqnarray*}
\textbf{Remark 3.} The discrete Jensen inequality has a long history, it has been applied in a huge number of works. In contrast, the discrete Wirtinger inequality for differences has been developed very recently: it has been published independently by several authors in slightly different forms in works \cite{hien15}, \cite{lee15}, \cite{Nam15a}, \cite{Nam15b}, \cite{seur15}; inequality (\ref{eg8uj})-(\ref{gy10uj}) is identical with that of  \cite{seur15}, and it is equivalent to all others.

\section{Stability analysis of discrete delayed systems}

The stability of the discrete delayed systems will be analyzed  in this section by a set of quadratic Lyapunov-Krasovskii functionals  (LKFs)
applying estimations
(\ref{eg25}). In the past decades, numerous different LKFs have been proposed and various techniques have been applied to reduce the conservatism of the results.
The effectiveness of the different LKFs are often tested by 
benchmark examples. Our purpose is
to evaluate systematically the reduction of conservativeness when applying estimations 
(\ref{eg25}). Similarly to
the useful approach of \cite{seur14b}, we will establish a hierarchy of linear matrix inequality (LMI) stability conditions based on
estimation of the forward difference of the LKF. These LMIs depend both on the parameter $m$ of the scalar product
(\ref{eg3}) and on the highest degree $\nu _m$ of the approximating orthogonal polynomials. The different cases corresponding to pairs
$(m,\nu_m)$ will be theoretically compared.

Consider the linear discrete time-delay system described by (\ref{uj9}).

Let us choose a positive integer $m$ and nonnegative integers $\nu_1 > \nu_2 > \ldots > \nu_{m} \geq 0 .$
Introduce the notation
\begin{equation}\label{gy70}
  x_{t-\tau}(i)=x(t- \tau +i), \hspace{0.5cm} \mbox{for} \hspace{0.5cm} i=0,1,\ldots,\tau-1,
\end{equation}
and set $\phi _j (t)= \sum _{i=0}^{\tau -1} p_{1j}(i)x_{t-\tau }(i).$
For $\nu _1 \geq 1,$ consider the extended and the augmented state variables
\begin{eqnarray}
\widetilde{x}(t)&= &\mbox{col}\left(x(t), \ \phi_0 (t), \ldots, \phi_{\nu_1-1} (t)\right) \in \mathbf{R}^{n_x(\nu _1+1)}, \nonumber\\                                              \widetilde{\Phi} (t)&= & \mbox{col}\left(x(t), \, x(t-\tau ), \, \frac{1}{\tau} \phi_0 (t), \ldots,  \frac{1}{\tau}  \phi_{\nu _1 -1} (t)
                             \right) \in \mathbf{R}^{n_x(\nu _1+2)},  \label{gy80}
\end{eqnarray}
respectively, and for $\nu_1 =0, $ set
\begin{eqnarray*}
\widetilde{x}(t) = x(t) , \hspace{0.4cm}   \widetilde{\Phi} (t)=  \mbox{col}\left(x(t), \, x(t-\tau )\right)  .
\end{eqnarray*}
Several further notations are needed for deriving the stability result.
Set
\begin{eqnarray*}
e_i &=&\left[0_{n_x\times n_x(i-1)}, \; I_{n_x}, \;  0_{n_x\times n_x(\nu _1+2-i)} \right], \; i=1,\ 2,\\ 
\mathcal{A} &=&  Ae_1+  A_d e_2, \hspace{4.0cm}
T_{\nu _1} = \mbox{diag}\left\{I_2, \; \tau  I_{\nu_1} \right\}, \\
\Gamma_{\nu_1} &=& \left(\mbox{diag}\left\{[1, \; 0], I_{\nu_1} \right\} T_{\nu _1}\right)\otimes I_{n_x}, \hspace{1.1cm}
\widetilde{\mathcal{Z}}_k = \mathcal{Z}_k T_{\nu _1},  \hspace{0.2cm} k=1,\ldots,m\\
\Lambda _l& = &\left(\left[ c_{1,l,1} \; c_{1,l,0} \; \lambda _{l,0} \; \ldots  \;  \lambda _{l,\nu_1 -1} \right]
T_{\nu _1}\right)\otimes I_{n_x},
\end{eqnarray*}
where  $\mathcal{Z}_k$ is defined by (\ref{gy50}),
$c_{1,l,0}=-p_{1l}(-1)$, $c_{1,l,1}=p_{1l}(\tau -1)$ as in the previous section, while $\lambda _{l,s}$ is defined by the
relation
\begin{equation*}
  p_{1l}(i-1)= \sum_{s=0}^{l}\lambda_{l,s}p_{1s}(i),
\end{equation*}
if $0 \leq s \leq l$ and $\lambda_{l,s}=0, $ if $l<s\leq\nu _1 .$
Finally set
\begin{eqnarray*}
\underline{\Lambda} _{\nu_1} = \left(
                     \begin{array}{llll}
                       \mathcal{A}^T &
                       \Lambda _{0}^T & \ldots &
                       \Lambda _{\nu_1 -1}^T
                     \end{array}
                   \right)^T .
\end{eqnarray*}

\begin{thm} \label{Th:3}
System (\ref{uj9}) is asymptotically stable, if there are 
matrices $P\in \mathbf{S}^{+}_{n_x(\nu_1 +1)}$, $Q, \ R_j \in \mathbf{S}^{+}_{n_x}$, $j=1,...,m$
 such that the LMI
\begin{eqnarray}
&&\Psi _{\nu _1}^{1}+ \Psi _{m,\nu _1}^{2} - \Psi _{m,\nu _1,...,\nu _m}^{3} < 0    \label{uj10},
\end{eqnarray}
has a solution, where
\begin{eqnarray}
&&\Psi _{\nu _1}^{1} = \underline{\Lambda}_{\nu _1}^{T} P \underline{\Lambda}_{\nu _1} - \Gamma_{\nu _1}^{T} P \Gamma_{\nu _1} + \mathcal{Q}_{\nu _1}   \label{uj10x}  \\
&&\Psi _{m,\nu _1}^{2} = \sum _{k=1}^{m} \left(  \begin{array}{c}
                                         \tau -1+k \\
                                         k
                                      \end{array}  \right) \left( \mathcal{A}-e_1  \right)^{T}R_k \left( \mathcal{A}-e_1  \right)  \label{uj10y} \\
&&\Psi _{m,\nu _1,...,\nu _m}^{3} = \sum _{k=1}^{m} \frac{1}{(k-1)!} \left( \widetilde{\mathcal{Z}}_k \otimes I \right)^{T} \mathcal{R}_k
\left( \widetilde{\mathcal{Z}}_k \otimes I \right)   \label{uj10z}
\end{eqnarray}
with
$\mathcal{Q}_{\nu _1}=\mbox{diag} \left\{ Q,\, -Q, 0_{\nu _1 \times \nu _1} \right\}  $ 
and $\mathcal{R}_k$ given by (\ref{gy50}).
\end{thm}

\textbf{Proof.}
Let us introduce
\begin{eqnarray*}
\rho _{t-\tau} (s)=x_{t+1-\tau}(s)-x_{t-\tau}(s)=x(t+1-\tau +s)-x(t-\tau +s),
\end{eqnarray*}
for $s=0,1,...,\tau -1$.
Consider the Lyapunov-Krasovskii functional candidate
\begin{eqnarray}
V(x_{t-\tau} , \rho _{t-\tau})=V_{00}(x_{t-\tau})+V_{01}(x_{t-\tau}) + \sum _{j=1}^{m}V_{1j}(\rho _{t-\tau}),        \label{plus1}
\end{eqnarray}
where
\begin{eqnarray}
V_{00}(x_{t-\tau})&=& \widetilde{x}^{T}(t)P\widetilde{x}(t),   \label{plus2}   \\
V_{01}(x_{t-\tau})&=& \sum _{s=0}^{\tau -1}x_{t- \tau}^{T}(s)Q x_{t- \tau} (s),  \label{plus3}   \\
V_{1,j}(\rho _{t-\tau})&=& \sum _{i_1 =0}^{\tau -1} ... \sum _{i_j =0}^{i_j -1} \sum _{s=i_j}^{\tau -1}
\rho _{t-\tau}^{T}(s)R_j \rho _{t-\tau} (s).    \label{uj9a}
\end{eqnarray}
Clearly, $V$ is positive definite.
Compute the forward difference of the LKF (\ref{plus1}) along the solution of (\ref{uj9}),
denoting the forward difference of $V$ by $\Delta v$ and the forward differences of terms
(\ref{plus2})-(\ref{uj9a}) by $\Delta v_{1j}.$
It can easily be verified that
\begin{eqnarray}
\widetilde{x}(t)&=&\Gamma_{\nu_1}  \widetilde{\Phi} (t).     \label{uj5} \\
\widetilde{x}(t+1)&=&\underline{\Lambda} _{\nu_1}  \widetilde{\Phi}(t),             \label{uj6}
\end{eqnarray}
therefore one obtains
\begin{eqnarray*}
\Delta v_{00}(t)&= & \widetilde{x}^{T}(t+1)P\widetilde{x}(t+1)-\widetilde{x}^{T}(t)P\widetilde{x}(t)=   \\
&& \widetilde{\Phi}^T (t) \left[ \underline{\Lambda} _{\nu_1}^{T} P \underline{\Lambda} _{\nu_1} - \Gamma _{\nu_1}^{T} P \Gamma _{\nu_1} \right]
\widetilde{\Phi}(t).
\end{eqnarray*}
Furthermore,
\begin{eqnarray*}
\Delta v_{01}(t)=x^{T}(t) Q x(t)-x^{T}(t-\tau) Q x(t-\tau)=\widetilde{\Phi}^T (t) \mathcal{Q}_{\nu_1} \widetilde{\Phi} (t),
\end{eqnarray*}
and
\begin{eqnarray*}
\Delta v_{1k}(t)&=& \sum _{i_1 =0}^{\tau-1} ... \sum _{i_k =0}^{i_k -1} \left( \rho _{t-\tau}^{T}(\tau ) R_k
\rho _{t-\tau} (\tau )- \rho _{t-\tau}^{T}(i_k) R_k \rho _{t-\tau} (i_k) \right) = \\
&=& \rho _{t-\tau}^{T}(\tau ) R_k \rho _{t-\tau}(\tau ) \sum _{i_1 =0}^{\tau-1} ... \sum _{i_k =0}^{i_k -1} 1-
\sum _{i_1 =0}^{\tau-1} ... \sum _{i_k =0}^{i_k -1} \rho _{t-\tau}^{T}(i_k ) R_k \rho _{t-\tau}(i_k ) =     \\
&=&   \left(
  \begin{array}{c}
    \tau - 1 +k \\
    k \\
  \end{array}
\right) \widetilde{\Phi}^{T} (t) \left( \mathcal{A}-e_1 \right)^{T} R_k \left( \mathcal{A}-e_1 \right)
\widetilde{\Phi} -J_k (\rho _{t-\tau}).
\end{eqnarray*}
For the estimation of $J_k (\rho _{t-\tau})$ we apply Theorem \ref{Th:2}. By taking into account (\ref{gy80}) and the definition of $\widetilde{\mathcal{Z}}_k,$ one obtains that
\begin{eqnarray*}
J_k (\rho _{t-\tau})\geq \frac{1}{(k-1)!} \widetilde{\Phi} ^T (t) \left( \widetilde{\mathcal{Z}} _k \otimes I \right)^T \mathcal{R}_k \left( \widetilde{\mathcal{Z}} _k \otimes I \right)
\widetilde{\Phi} (t).
\end{eqnarray*}
Therefore,
\begin{eqnarray*}
\Delta v(t) \leq \widetilde{\Phi} ^T (t) \left\{ \Psi _{\nu _1}^{1} + \Psi _{m,\nu _1}^{2} -
\Psi _{m,\ \nu _1,...,\nu _m}^{3} \right\} \widetilde{\Phi}(t),
\end{eqnarray*}
i.e. if (\ref{uj10}) holds true, there exist a constant $\varepsilon >0$ such that $\Delta v(t) \leq -\varepsilon\|x(t)\|^2,$ thus the statement of the theorem follows.
$\Box$

\textbf{Remark 4.} Matrices $\underline{\Lambda} _{\nu_1}$, $\mathcal{Z} _k,$ together with the generating  Wolfram Mathematica programs are given in \cite{kisnagy} for several $\nu_1, $ $k$ and $\nu_k$.

\section{Hierarchy of the LMI stability conditions}

This section is devoted to proving that the stability conditions (\ref{uj10}) obtained for different choices of $m$ and
$\nu_1 > \nu_2 > \ldots > \nu_{m} \geq 0 $ can be arranged into a hierarchy table.

\begin{defi}
Let $\mathcal{L}^{(i)}$ denote LMI (\ref{uj10}) obtained with the choice of $m^{(i)}, \nu_1^{(i)},  \nu_2^{(i)},  \ldots,  \nu_{m^{(i)}}^{(i)}, $ ($i=1,2$).  \emph{LMI $\mathcal{L}^{(2)}$ outperform  LMI $\mathcal{L}^{(1)}$} 
 if and only if for any $\tau $, for which $\mathcal{L}^{(1)}$
has a feasible solution, LMI $\mathcal{L}^{(2)}$ has it, as well. This relation will be denoted by $\mathcal{L}^{(1)} \prec \mathcal{L}^{(2)}.$ 
\end{defi}

Let symbol $\mathcal{L}_{\nu_1,...,\nu _m}^{m}$ denote LMI (\ref{uj10}). Let us arrange these LMIs into a hierarchy table as follows:
\begin{equation*}
  \begin{array}{lllll}
  \mathcal{L}_0^1 & \mathcal{L}_1^1 & \mathcal{L}_2^1 & ... & \mathcal{L}_{\nu_1}^1 \\
    & \mathcal{L}_{1,0}^2 & \mathcal{L}_{2,1}^2 & ... & \mathcal{L}_{\nu_1,\nu_1-1 }^2 \\
    &   & \mathcal{L}_{2,1,0}^3 & ... & \mathcal{L}_{\nu_1,\nu_1-1,\nu_1-2 }^3 \\
   &  &  & \ddots &  \\
   &  &  &  & \mathcal{L}_{\nu_1,\nu_1-1,...,\nu_1-(m-1) }^m
\end{array}
\end{equation*}

 Note that this table is finite for any given $\tau$ both to the right and to the bottom, since $\tau$ determines the maximal degree of the orthogonal polynomials that can be considered.

We want to show that a certain LMI in this table outperforms any other LMI, which is situated above and/or to the left of it.

\begin{thm} \label{Th:4}
Let the integers $1 \leq m$ and $m < \nu_1^{\ast}$ be given. Then
\begin{itemize}
\item
for $1 \leq \ell \leq m$ and $\ell-1 \leq \nu_1 < \nu_1^{\ast},$
\begin{equation}\label{gy100}
  \mathcal{L}_{\nu_1,\nu_1-1,...,\nu_1+1-\ell }^\ell  \prec \mathcal{L}_{\nu_1+1,\nu_1,...,\nu_1+2-\ell }^\ell \ ;
  \end{equation}
\item
for $1 \leq \ell < m$ and $\ell-1 \leq \nu_1 \leq \nu_1^{\ast},$
 \begin{equation}\label{gy110}
  \mathcal{L}_{\nu_1,\nu_1-1,...,\nu_1+1-\ell }^\ell  \prec \mathcal{L}_{\nu_1,\nu_1-1,...,\nu_1+1-\ell,\nu_1-\ell}^{\ell+1} \ .
\end{equation}
\end{itemize}
\end{thm}

\textbf{Proof.}
\emph{Part I.}
First we show that (\ref{gy100}) is valid. Suppose that $\mathcal{L}_{\nu_1,\nu_1-1,...,\nu_1+1-\ell }^\ell$ has the feasible solution
$P^{(1)},$ $Q^{(1)}, $ $R_1^{(1)}, \ldots, R_\ell^{(1)} $ for a given $\tau.$
 We show that, for the same value of $\tau,$ there exists a feasible solution $P^{(2)},$ $Q^{(2)}, $ $R_1^{(2)}, \ldots, R_\ell^{(2)} $ of
 $\mathcal{L}_{\nu_1+1,\nu_1,...,\nu_1+2-\ell }^\ell.$ In what follows, we shall use the upper index "$(1)$" and "$(2)$" to distinguish the matrices occurring in the first and in the second LMI, respectively.

Let us  seek $P^{(2)}$ in the form of $P^{(2)}=\mbox{diag}\left\{ P^{(1)}, \ \varepsilon I_{n_x}\right\},$ where $\varepsilon$ is some positive constant, while we set $Q^{(2)}=Q^{(1)}, $ $R_1^{(2)}=R_1^{(1)}, \ldots, R_\ell^{(2)}=R_\ell^{(1)}. $

Calculate first matrix $\Psi_{\nu_1+1}^{1}.$  Observe that
\begin{eqnarray}
\mathcal{Q}_{\nu_1 +1} = \begin{bmatrix}
      \mathcal{Q}_{\nu_1} & 0 \\
                   0 & 0 \\
  \end{bmatrix}, \hspace{0.3cm} \Gamma_{\nu_1 +1}=\begin{bmatrix}
      \Gamma_{\nu_1} & 0 \\
                   0 & \tau I_{n_x} \\
  \end{bmatrix}, \hspace{0.3cm} \underline{\Lambda} _{\nu_1 +1}=\begin{bmatrix}
      \underline{\Lambda} _{\nu_1} & 0 \\
         \Lambda_{\nu_1 +1,1} & \Lambda _{\nu_1 +1,2} \
  \end{bmatrix},
\end{eqnarray}
where
\begin{equation*}
    \Lambda_{\nu_1 +1,1}=\left( c_{1,\nu_1 ,1}, \; \; c_{1,\nu_1 ,0}, \; \; \lambda_{\nu_1 ,0}, \; \;  \ldots \; ,\lambda_{\nu_1 -1} \right)T_{\nu1} \otimes I  ,\hspace{0.5cm} \mbox{and} \hspace{0.5cm}
    \Lambda_{\nu_1 +1,2} = \tau I_{n_x}.
\end{equation*}
By straightforward calculation we obtain
\begin{eqnarray}
&&\underline{\Lambda}_{\nu_1 +1}^{T}P^{(2)} \underline{\Lambda}_{\nu_1 +1} =
\begin{bmatrix}
      \underline{\Lambda}_{\nu_1}^{T}P^{(1)} \underline{\Lambda}_{\nu_1} & 0 \\
                   0 & 0 \\
  \end{bmatrix} +
  \varepsilon \begin{bmatrix}
      0 & \left( \Lambda_{\nu_1 +1}^{(1)} \right)^{T}  \\
                   0 & (\Lambda ^{(2)})^{T} \\
  \end{bmatrix}
  \begin{bmatrix}
      0 & 0  \\
                   \Lambda_{\nu_1 +1}^{(1)} & \Lambda ^{(2)} \\
  \end{bmatrix}
    \label{gy210}    \\
&&\Gamma_{\nu_1 +1}^{T} P^{(2)} \Gamma_{\nu_1 +1} =
\begin{bmatrix}
     \Gamma_{\nu_1}^{T} P^{(1)} \Gamma_{\nu_1} & 0  \\
                  0 & \varepsilon \tau I \\
  \end{bmatrix} \label{gy220}
%
\end{eqnarray}

As far as $\Psi _{\ell,\nu _1+1}^{2}$ is concerned observe that
\begin{eqnarray}
\Psi _{\ell,\nu _1+1}^{2} = \begin{bmatrix}
     \Psi _{\ell,\nu _1}^{2} & 0 \\
                             0 & 0 \\
  \end{bmatrix}.
\end{eqnarray}
Finally consider matrix $\Psi _{\ell,\nu_1+1,\nu_1,...,\nu_1+2-\ell }^{3}:$ 
\begin{eqnarray*}
\Psi _{\ell,\nu _1+1,\nu_1,...,\nu_1+2-\ell}^{3}= \sum _{k=1}^{\ell} \frac{1}{(k-1)!}
\left( \widetilde{\mathcal{Z}}_k^{(2)} \otimes I \right)^{T} \mathcal{R}_k^{(2)}
\left( \widetilde{\mathcal{Z}}_k ^{(2)}\otimes I \right),
\end{eqnarray*}
where
\begin{equation*}
\mathcal{R}_k ^{(2)} = \begin{bmatrix}
\mathcal{R}_k ^{(1)} & 0 \\
       0    & \chi_{k,\nu_1-k+2}R_k \\
  \end{bmatrix}, \hspace{1cm}
  \widetilde{\mathcal{Z}}_k ^{(2)}=\begin{bmatrix}
     \widetilde{\mathcal{Z}}_k ^{(1)} & 0 \\
       \widetilde{\underline{\zeta}}_{\nu_1-k+2,1}^{k}    & \widetilde{{\zeta}}_{\nu_1-k+2,2}^{k} \\
  \end{bmatrix}
\end{equation*}
with
\begin{eqnarray*}
 && \widetilde{\underline{\zeta}}_{\nu_1-k+2,1}^{k}=
  \left( c_{k,\nu_1-k+2,1},  \; \; c_{k,\nu_1-k+2,0},  \; \; \zeta_{\nu_1-k+2,0}^{k}, \;\; \ldots , \;\;
  \zeta_{\nu_1-k+2,\nu_1-1}^{k} \right)T_{\nu1} \\&&
  \hspace*{-2.1cm}\mbox{and}\\
  && \widetilde{{\zeta}}_{\nu_1-k+2,2}^{k}=\zeta_{\nu_1-k+2,\nu_1}^{k} \tau.
\end{eqnarray*}
We recall that $\zeta_{\nu_1-k+2,\nu_1}^{k}\neq 0,$ (see the comment above (\ref{gy50})).
By straightforward calculation we obtain
\begin{eqnarray}
&&(\widetilde{\mathcal{Z}}_k ^{(2)} \otimes I)^{T} \mathcal{R}_k ^{(2)} (\widetilde{\mathcal{Z}}_k ^{(2)} \otimes I) =
\begin{bmatrix}
(\widetilde{\mathcal{Z}}_k^{(1)}  \otimes I)^{T} \mathcal{R}_k ^{(1)}(\widetilde{\mathcal{Z}}_k ^{(1)} \otimes I) & 0 \\
    0 & 0 \\
  \end{bmatrix}
  \nonumber \\
&& \hspace{1.5cm}+
\begin{bmatrix}
       (\widetilde{\underline{\zeta}}_{\nu_1-k+2,1}^{k} \otimes I)^{T} \\
     \widetilde{{\zeta}}_{\nu_1-k+2,2}^{k} I \\
\end{bmatrix}
             \chi_{k,\nu_1-k+2}R_k
\begin{bmatrix}
                 \widetilde{\underline{\zeta}}_{\nu_1-k+2,1}^{k} \otimes I & \widetilde{{\zeta}}_{\nu_1-k+2,2}^{k} I \\
\end{bmatrix} . \label{gy200}
\end{eqnarray}
Since $\chi_{k,\nu_1-k+2}R_k>0,$  the second term on the right hand side of (\ref{gy200}) is nonnegative. Therefore,
 matrix $\Psi _{\ell,\nu_1+1,\nu_1,...,\nu_1+2-\ell }^{3}$ can be estimated from below by keeping only one of the second terms at summation,
e.g. when $k=1$. One obtains
\begin{eqnarray}
 \Psi _{\ell,\nu_1+1,\nu_1,...,\nu_1+2-\ell }^{3} &\geq&
  \begin{bmatrix}
    \Psi _{\ell,\nu_1,\nu_1-1,...,\nu_1+1-\ell }^{3} & 0 \\
    0 & 0 \\
  \end{bmatrix}
  +    \nonumber   \\
&&
\begin{bmatrix}
       (\widetilde{\underline{\zeta}}_{\nu_1+1,1}^{1} \otimes I)^{T} \\
     \widetilde{{\zeta}}_{\nu_1+1,2}^{1} I \\
\end{bmatrix}
             \chi_{1,\nu_1+1}R_1
\begin{bmatrix}
                 \widetilde{\underline{\zeta}}_{\nu_1+1,1}^{1} \otimes I & \widetilde{{\zeta}}_{\nu_1+1,2}^{1} I \\
\end{bmatrix} .                          \label{uj14}
\end{eqnarray}
To be short, let us use the notations $$\mathcal{M}^{(1)}=\Psi _{\nu _1}^{1}+ \Psi _{\ell, \nu _1}^{2} -
\Psi _{\ell,\nu_1,\nu_1-1,...,\nu_1+1-\ell }^{3}$$ and $$\mathcal{M}^{(2)}=\Psi _{\nu _1+1}^{1}+ \Psi _{\ell, \nu _1+1}^{2} -
\Psi _{\ell,\nu_1+1,\nu_1,...,\nu_1+2-\ell}^{3}.$$
Employing (\ref{gy210})-(\ref{uj14}) we obtain for $\mathcal{L}_{\nu_m +m,...,\nu_m +2, \nu_m +1}$ that
\begin{eqnarray}
\mathcal{M}^{(2)} &\leq &
 \begin{bmatrix}
    I & (\widetilde{\underline{\zeta}}_{\nu_1+1,1}^{1} \otimes I)^{T} \\
    0 & \widetilde{{\zeta}}_{\nu_1+1,2}^{1} I\\
  \end{bmatrix}
  \begin{bmatrix}
    \mathcal{M}^{(1)} & 0 \\
    0 & -\chi_{1,\nu_1+1}R_1  \\
  \end{bmatrix}
 \begin{bmatrix}
    I & 0 \\
    \widetilde{\underline{\zeta}}_{\nu_1+1,1}^{1} \otimes I & \widetilde{{\zeta}}_{\nu_1+1,2}^{1} I\\
  \end{bmatrix}   \nonumber    \\
&& \hspace{1.7cm} +\varepsilon \begin{bmatrix}
      0 & \left( \Lambda_{\nu_1 +1}^{(1)} \right)^{T}  \\
                   0 & (\Lambda ^{(2)})^{T} \\
  \end{bmatrix}
  \begin{bmatrix}
      0 & 0  \\
                   \Lambda_{\nu_1 +1}^{(1)} & \Lambda ^{(2)} \\
  \end{bmatrix} .     \label{uj15}
\end{eqnarray}
Since $\widetilde{{\zeta}}_{\nu_1+1,2}^{1}\neq 0,$ the two extreme multiplier matrices in the first term on the right hand side of (\ref{uj15}) are invertible. Therefore, if $\mathcal{L}_{\nu_1,\nu_1-1,...,\nu_1+1-\ell }^\ell$ has a feasible solution, there exists
a $\mu_1 >0$ such that the first term on the right hand side of (\ref{uj15}) is less
than $-\mu_1 I.$
On the other hand, there exists a constant $\lambda^{*}$ such that the matrix product on the right hand side of (\ref{uj15}) can be estimated as
\begin{eqnarray*}
\begin{bmatrix}
      0 & \left( \Lambda_{\nu_1 +1}^{(1)} \right)^{T}  \\
                   0 & (\Lambda ^{(2)})^{T} \\
  \end{bmatrix}
  \begin{bmatrix}
      0 & 0  \\
                   \Lambda_{\nu_1 +1}^{(1)} & \Lambda ^{(2)} \\
  \end{bmatrix} \leq \lambda^{*} I,
\end{eqnarray*}
 therefore $\mathcal{M}^{(2)}<0$ is satisfied, if $\varepsilon \lambda^{*} < \mu_1 .$ This proofs that (\ref{gy100}) is true.

\emph{Part II.}
Next we show that one can move  downwards in the hierarchy, too, if $\ell$ is increased with fixed $\nu_1$,
i.e. $\mathcal{L}_{\nu_1,\nu_1-1,...,\nu_1+1-\ell }^\ell  \prec \mathcal{L}_{\nu_1,\nu_1-1,...,\nu_1+1-\ell,\nu_1-\ell}^{\ell+1} .$
We recalculate the terms (\ref{uj10x})-(\ref{uj10z}) again in this case. We obtain
\begin{eqnarray*}
&&\Psi_{\ell+1,\nu_1}^{1} = \Psi_{\ell,\nu_1}^{1},   \\
&&\Psi_{\ell+1,\nu_1}^{2} = \Psi_{\ell,\nu_1}^{2} + \left(
                                              \begin{array}{c}
                                                \tau -1+\ell+1 \\
                                                \ell+1 \\
                                              \end{array}
                                            \right) (\mathcal{A}-e_1)^{T} R_{\ell+1}(\mathcal{A}-e_1),    \\
&&\Psi_{\ell+1,\nu_1,\nu_1-1,...,\nu_1+1-\ell,\nu_1-\ell}^{3} = \Psi_{\ell,\nu_1,\nu_1-1,...,\nu_1+1-\ell  }^{3} +
\frac{1}{\ell!} \left(\widetilde{\mathcal{Z}}_{\ell+1} \otimes I \right)^T \mathcal{R}_{\ell+1} \left( \widetilde{\mathcal{Z}}_{\ell+1} \otimes I \right) .
\end{eqnarray*}
Analogously to Part I, we shall use the brief notation
$$\mathcal{N}^{(1)}=
\Psi_{\ell,\nu_1}^{1} + \Psi_{\ell+1,\nu_1}^{2} -
\Psi_{\ell,\nu_1,\nu_1-1,...,\nu_1+1-\ell  }^{3}$$
and
$$\mathcal{N}^{(2)}=\Psi_{\ell+1,\nu_1}^{1}+ \Psi _{\ell+1, \nu _1}^{2} -
\Psi _{\ell+1,\nu_1,\nu_1-1,...,\nu_1+1-\ell,\nu_1-\ell }^{3}.$$
Then $\mathcal{N}^{(2)}$ can be expressed as $\mathcal{N}^{(2)}= \mathcal{N}^{(1)}+\Upsilon, $ where
\begin{eqnarray*}
\Upsilon& = &   \left(
                                              \begin{array}{c}
                                                \tau -1+\ell+1 \\
                                                \ell+1 \\
                                              \end{array}
                                            \right) (\mathcal{A}-e_1)^{T} R_{m+1}(\mathcal{A}-e_1)    \\
&& - \; \frac{1}{\ell!} \left(\widetilde{\mathcal{Z}}_{\ell+1} \otimes I \right)^T \mathcal{R}_{\ell+1} \left( \widetilde{\mathcal{Z}}_{\ell+1} \otimes I \right).                                            \end{eqnarray*}
Let us seek $R_{\ell+1}$ in the form of $R_{\ell+1}=\varepsilon I.$ Then there is a $\lambda^{**}$ such that $\Upsilon$ can be estimated as
\begin{eqnarray*}
\Upsilon &=& \varepsilon \left\{\left(   \begin{array}{c}
                                                \tau -1+m+1 \\
                                                m+1 \\
                                              \end{array}
                                            \right) (\mathcal{A}-e_1)^{T} (\mathcal{A}-e_1)  \right.   \\
&&\left. - \frac{1}{\ell!} \left(\widetilde{\mathcal{Z}}_{\ell+1} \otimes I \right)^T
\mbox{diag} \left\{ \chi _{\ell,0}I, \; \ldots , \;
\chi _{\ell+1,\nu_1+1-\ell}I  \right\}
\left( \widetilde{\mathcal{Z}}_{\ell+1} \otimes I \right) \right\}\leq \varepsilon \lambda^{**}I.
\end{eqnarray*}
If $\mathcal{L}_{\nu_1,\nu_1-1,...,\nu_1+1-\ell }^\ell$ has a feasible solution, then there exists a $\mu _2>0$ such that $$\mathcal{N}^{(1)}<-\mu _2 I,$$
therefore  $\mathcal{N}^{(2)}<0$
is satisfied if $ \varepsilon \lambda^{**} < \mu _2 .$  This proofs that (\ref{gy110}) is valid.
$\Box$

\section{Numerical examples}

In this section, we apply the proposed method to two benchmark examples that have been extensively used in the literature to compare the results. The third example is a slight modification of Example 1 investigated in \cite{gyt15} in a different situation.

\subsection{Some remarks on the implementation}
The computations of the data matrices for the application of Theorem \ref{Th:3}
have been performed by using Wolfram Mathematica. First the monic orthogonal polynomials have been computed, then the polynomials of indeces $(k,j)$ have been normalized with a multiplier $\pi_j$ so that $p_{1j}(-1)=(-1)^j$ is satisfied. In this way, it has been achieved that the elements of all matrices have values of reasonable order of magnitude, and the computations  are 
numerically stable. As a result, we have e.g. the norm-squares $\|p_{1,0}\|^2=\tau$ and  $\|p_{1,j}\|^2=\frac{\tau}{2j+1}\prod _{i=1}^{j}\frac{\tau-i}{\tau+i},$ if $j=1,...,\nu_1, $ while
matrix $\Lambda_5$ is as follows:
\begin{equation*}
  \Lambda_5=\begin{bmatrix}
      1 & -1 & 1 & 0 & 0 & 0 & 0 \\
      \frac{\tau-1}{\tau+1} & 1 & \frac{-2}{\tau+1} & 1 & 0 & 0 & 0 \\
      \prod _{i=1}^{2}\frac{\tau-i}{\tau+i} & -1 & \frac{6}{\prod _{i=1}^{2}(\tau+i)} & \frac{-6}{\tau+2} & 1 & 0 & 0 \\
      \prod _{i=1}^{3}\frac{\tau-i}{\tau+i} & 1 & \frac{-2(\tau^2+11}{\prod _{i=1}^{3}(\tau+i)} & \frac{30}{\prod _{i=2}^{3}(\tau+i)}
      & \frac{-10}{\tau+3} & 1 & 0 \\
      \prod _{i=1}^{4}\frac{\tau-i}{\tau+i} & -1 & \frac{20(\tau^2+5)}{\prod _{i=1}^{4}(\tau+i)} &
      \frac{-6(\tau^2+26)}{\prod _{i=2}^{4}(\tau+i)} & \frac{70}{\prod _{i=3}^{4}(\tau+i)} & \frac{-14}{\tau+4} & 1 \\
      \prod _{i=1}^{5}\frac{\tau-i}{\tau+i} & 1 & \frac{-2(\tau^4+85 \tau^2 + 274)}{\prod _{i=1}^{5}(\tau+i)} &
      \frac{84(\tau^2+11)}{\prod _{i=2}^{5}(\tau+i)} & \frac{-10(\tau^2 + 47)}{\prod _{i=3}^{5}(\tau+i)} & \frac{126}{\prod _{i=4}^{5}(\tau+i)} &
      \frac{-18}{\tau +5} \\
  \end{bmatrix}
\end{equation*}
%
The norm-squares of the other polynomials, the matrices $\mathcal{Z}_k$ as well as the details of computations can be found in \cite{kisnagy}.

The LMIs have been solved by using MATLAB LMI Toolbox. The computations have been performed for $m=1,...,4$ and $\nu_1=0,...,5,$ but tables
  below include only the most informative results. 

\subsection{Numerical experiments}
\textbf{Example 1.} 
Consider system (\ref{uj9}) with
\begin{equation*}
  A= \begin{bmatrix}
      0.8 & 0  \\
        0 & 0.91 \\
  \end{bmatrix}, \hspace{1cm}
  A_d= \begin{bmatrix}
      -0.1 & 0  \\
        -0.1 & -0.1 \\
  \end{bmatrix}.
\end{equation*}
The analytical range of the delay that retains stability of the system is $[0, 58]\cap \mathbf{Z}.$
The number of decision variables is $7021$, if the discrete Lyapunov inequality is used to determine the analytical bound.
The results obtained by methods proposed in one of the most recent references \cite{ZhH} and by Theorem \ref{Th:3} for different values of $m$ and $\nu _1$  are given in Table 1. Further comparisons with results of several recent references is given in \cite{ZhH}. For $m>1,$  the values of $ \nu_j$ are set as $\nu_j=\nu_1 - (j-1).$

   \begin{table}[!ht] \label{Tab:1}
\caption{{\footnotesize Delay upper bound for Example 1 and Example 2 }}
{\small
\begin{center}
\begin{tabular}{lcccccccccccc}
 \toprule
Example & \hspace{0.5cm} & & &  &1 & & &\hspace{0.3cm} & &  & 2& \\
\hline
Method & & & $m$ & $\nu_1$ & $\tau_M$ & NoDVs & & & $m$ & $\nu_1$ & $\tau_M$ & NoDVs\\
\hline
 Zhang et al. \cite{ZhH}  & &  & 1 & 1 & 57  &$16$& & & $1$ & $1$ & $151$  &$16$ \\
Nam et al.  \cite{Nam15b}  & &  &  &  &   & & &  & $2$&  $2$ & $168$  & $ ^{(*)}$ \\
  Theorem \ref{Th:3} & & & 1  & 0 & 42     & $9$& & & $1$  & $1$ & $151$    & $16$ \\
  & & & 1  & 1 & 57     & $16$                & &  & $1$  & $2$ & $168$     & $27$\\
   & &  & 2  & 1 & 57     & $19$              & &   & $2$  & $2$ & $168$     & $30$\\
   & & & 1  & 2 & 58     & $27$               & &  & $1$  & $4$ & $169$     & $61$    \\
   & &  & 2  & 2 & 58     & $30$              & &  & $2$  & $4$ & $169$     & $64$ \\
 \bottomrule
        {\footnotesize $ ^{(*)}$ not available} & & & & & & & & & &  \\
 \end{tabular}  \\
 \end{center}
  \vskip-2mm
     }
   \end{table}

\textbf{Example 2.} 
Consider system (\ref{uj9}) with
\begin{equation*}
  A= \begin{bmatrix}
      1 & 0.01  \\
        -0.02 & 1.001 \\
  \end{bmatrix}, \hspace{1cm}
  A_d= \begin{bmatrix}
      0 & 0  \\
         0.01 & 0 \\
  \end{bmatrix}.
\end{equation*}
The analytical range of the delay that retains stability of the system is $[12, 169]\cap \mathbf{Z}.$
The number of decision variables is $57970$, if the discrete Lyapunov inequality is used to determine the analytical bound.
The results obtained by methods proposed in most recent references \cite{ZhH} and \cite{Nam15b} and by Theorem \ref{Th:3} for different values of $m$ and $\nu _1$  are given in Table 1.
 Similarly to the previous example, for $m>1,$  the values of $ \nu_j$ are set as $\nu_j=\nu_1 - (j-1).$ Further comparisons with  results of several recent references is given in \cite{Nam15b}.
  The LMIs of Theorem \ref{Th:3} were feasible for $\tau=12$ in all cases.

 \textbf{Example 3.} 
 Consider system (\ref{uj9}) with
\begin{equation*}
  A= \begin{bmatrix}
      0.12 & 0 & -0.12  \\
      0.06 & 0.36 & 0 \\
      0 & 0.24 & 0.72 \\
  \end{bmatrix}, \hspace{1cm}
  A_d= \begin{bmatrix}
      -0.4 & 0 & 0  \\
      0 & -0.2 & 0.2 \\
      0 & 0 & -0.4 \\
  \end{bmatrix}.
\end{equation*}
The analytical range of the delay that retains stability of the system is $[0, 56]\cap \mathbf{Z}.$
The number of decision variables is $14706$, if the discrete Lyapunov inequality is used to determine the analytical bound.
Similarly to the previous example, for $m>1,$  the values of $ \nu_j$ are set as $\nu_j=\nu_1 - (j-1).$
The results obtained by methods proposed in \cite{ZhH} and by Theorem \ref{Th:3} for $m=1$ and for different values
 of $\nu _1$ are given in Table 2. In cases of $m=2,3,4,5$, the delay bounds were found to be the same as for $m=1$ applying
 the same value of $\nu_1.$

\begin{table}[!ht] \label{Tab:2}
\caption{{\footnotesize Delay upper bound for Example 3}}
{\small
\begin{center}
\begin{tabular}{lcccccc}
\toprule
Method  & \hspace{0.3cm} & & $m$ & $\nu_1$ & $\tau_M$ & NoDVs \\
\hline
 Zhang et al. 2015 \cite{ZhH} & &  & $1$ & $1$ & $50$  &$33$ \\
   \hline
     Theorem \ref{Th:3} & & & $1$  & $0$ & $34$    & 18 \\
  & & & $1$  & $1$ & $50$     & 33 \\
   & &  & $1$  & $2$ & $52$     & 57 \\
  & & & $1$  & $3$ & $52$     & 90 \\
  & &  & $1$  & $4$ & $55$     & 132 \\
   & & & $1$  & $5$ & $56$     & 183 \\
   \bottomrule
 \end{tabular}  \\
 \end{center}
  \vskip-2mm
     }
   \end{table}

\subsection{Discussion}
The numerical examples show that the application of the discrete Wirtinger inequality (case $m=1$ and $\nu_1=1$) reduces the conservativeness of results obtained by the application of Jensen's inequality (case $m=1$ and $\nu_1=0$). They also show that further improvement can be achieved by the application of the inequalities of Theorem \ref{Th:2} both via increasing the number of multiple summation terms in the LKF (case $m>1$)  and the improvement of the lower estimations via increasing the
 degree of the orthogonal polynomials (i.e. $\nu_1>1 , \,..., \nu_m \geq 0 $). 
Beside the above examples, we tested our approach for several other examples from the literature with the experience
as follows. The increase of the complexity of the LKF (i.e. the increase of $m$) did not resulted in a better delay bound than the LKF with
$m=1,$ if the same $\nu _1 $ was applied. This means that the improvement is primarily due to the increase of the dimension of the extended
state variable. This does not contradict to the reported improvements in the case the application of triple, etc. summation terms in the LKF,
since  - on the one hand - the applied lower estimations lead to introduction of some extended state variables with increased dimension.
On the other hand, several authors
apply in their developments not only a more complex LKF, but other methods as well (e.g. adding 'zero equality', see e.g.
\cite{Nam15b}, relaxation
of the requirement of positive definiteness of certain matrices in the LKF, see e.g. \cite{XuDiscr13}, etc.).
We have not applied these latter methods, since we wanted only to investigate the effect of the improvement of the lower estimation and the effect of the application of multiple summation terms in the LKF.
In sum, the increase of the dimension of the extended state variable plays the basic role in the improvement.

We emphasize that analytical delay bounds could be achieved in all examples, if a sufficiently tight lower estimation is applied. Apparently, the necessary number of decision variables is dramatically lower than that under the application of the necessary and
sufficient condition of stability.

\section{Conclusion}

In this paper, multiple summation inequalities were presented in the case of arbitrary number of summation both for functions and differences.
The new inequalities involve the discrete Jensen's and Wirtinger's inequalities, as well as the recently presented inequalities for single and double summation in \cite{Nam15b}. Applying the obtained inequalities, a set of sufficient LMI stability conditions for linear discrete-time delay systems are derived. It was proven that these LMI conditions could be arranged into a bidirectional hierarchy establishing a rigorous theoretical basis for comparison of conservatism of the investigated methods.
It was shown by some benchmark examples that the proposed methods give better upper bounds for the tolerable time delay than the best ones that we could find in the previously published literature. Several numerical examples showed that the improvement of the lower estimations
can result in as much improvements as the application of more complex LKFs.



\begin{thebibliography}{99}

\bibitem
{Briat14} C. Briat,  \emph{ Linear parameter-varying and time-delay systems.} Springer, 2014.

\bibitem
{frid2014} E. Fridman, \emph{Introduction to time-delay systems: Analysis and control.} Springer, 2014.

\bibitem
{gau}
W. Gautschi,
Orthogonal polynomials computation and approximation, Oxford University
Press, 2004.

\bibitem
{gye15} \'{E}. Gyurkovics,
A note on Wirtinger-type integral inequalities for time-delay systems,
Automatica 61 November (2015) 44-46.

\bibitem
{gyt15} \'{E}. Gyurkovics, T. Tak\'acs,
Robust dynamic output feedback guaranteed cost control for discrete-time systems with time-varying delays,
Asian Journal of Control 17 2 (2015) 687-698.

\bibitem{Hardy} G.L. Hardy, J.E. Littlewood, G. Polya,
Inequalities,
Cambridge, UK: Cambridge University Press, (1934).

\bibitem
{hien15} L.V. Hien, H. Trinh,
New summation inequalities and their applications to discrete-time delay systems,
{\tt arXiv:1505.06344v1 [math.OC]}

\bibitem{kisnagy} K. Kiss, I. Nagy, Computations of discrete ortogonal polynomials with Wolfram Mathematica,
   http://www.math.bme.hu/\verb"~"kk/Mathematica\verb"_"notebooks\verb"_"KK.htm 
   (Available on 05.12.2015.)



\bibitem{marek} M. Kuczma, An Introduction to the Theory of Functional Equations and Inequalities.
Cauchy's Equation and Jensen's Inequality, Birkhäuser Verlag AG, Basel – Boston – Berlin
(2nd edition) (2009).

\bibitem
{kwon14}
 O.M. Kwon, M.J. Park, J.H. Park, S.M. Lee, E.J. Cha,
On less conservative stability criteria for neural networks with time-varying
delays utilizing Wirtinger-based integral inequality,
Hindawi Publishing Corporation, Mathematical Problems in Engineering,
2014 ID ID 859736 http://dx.doi.org/10.1155/2014/859736 (2014)

\bibitem
{lee15} W.I. Lee, P.G. Park, S.Y. Lee, R.W. Newcomb,
Auxiliary function-based summation inequalities for quadratic functions and
their application to discrete-time delay systems, IFAC-PapersOnLine
48-12 (2015) 203-208.

\bibitem{LFrAut12} K. Liu,  E. Fridman,  Wirtinger's inequality and Lyapunov-based sampled-data
stabilization. Automatica,  48 (2012) 102-108.

\bibitem
{liu15} X. Liu, X. Zhao,
Stability analysis of discrete-time switched systems: a switched homogeneous
Lyapunov function method, to appear in International Journal of Control, (2015)
DOI: 10.1080/00207179.2015.1075254.

\bibitem
{Nam15a} P.T. Nam, P.N. Pathirana, H. Trinh,
Discrete Wirtinger-based inequality and its application,
Journal of the Franklin Institute,
352.5 (2015) 1893-1905.

\bibitem
{Nam15b} P.T. Nam, H. Trinh, P.N. Pathirana,
Discrete inequalities based on a multiple auxiliary functions and their
applications to stability analysis of time-delay systems,
Journal of the Franklin Institute, (2015)
http://dx.doi.org/10.1016/j.jfranklin.2015.09.018.

\bibitem
{park15} M. Park, O.M. Kwon, J.H. Park, S.M. Lee, E. Cha,
Stability of time-delay systems via Wirtinger-based integral inequality,
Automatica 55 (2015) 204-208.


\bibitem{KwL}
P. Park,  Won Il Lee,  S.Y. Lee, Auxiliary function-based integral inequalities for quadratic functions and their applications to time-delay systems. Journal of the Franklin Institute 352.4 (2015): 1378-1396.

\bibitem
{seur12} A. Seuret, F. Gouaisbaut,
On the use of the Wirtinger inequalities for time-delay
systems, 10th IFAC Workshop on Time Delay Systems, June 2012, Boston,
United States, (2012)
https://hal.archives-ouvertes.fr/hal-00697498, 15 May.

\bibitem
{seur14} A. Seuret, F. Gouaisbaut,
Complete quadratic Lyapunov functionals using Bessel-Legendre inequality,
European Control Conference, June 2014, Strasbourg, France, (2014)
https://hal.archives-ouvertes.fr/hal-00983321v2, 22 September.

\bibitem
{seur14b} A. Seuret, F. Gouaisbaut,
Hierarchy of LMI conditions for the stability analysis of time delay
systems, Systems \& Control Letters 81 (2015) 1-7.

\bibitem
{seur15} A. Seuret, F. Gouaisbaut, E. Fridman,
Stability of discrete-time systems with time-varying delays via a novel
summation inequality, IEEE Transactions on Automatic Control, 60 (2015) 2740 - 2745.

\bibitem
{XuDiscr13} S. Xu, J. Lam, B. Zhang, Y. Zou, A new result on the delay-dependent stability of discrete systems with time-varying delays. International Journal of Robust and Nonlinear Control  24 (2014) 2512-2521.

\bibitem{WHSh} M. Wu M, Y. He, J-H. She, Stability Analysis and Robust Control of Time-Delay Systems,
Science Press Beijing and Springer-Verlag Berlin Heidelberg, 2010.


\bibitem
{ze15} H.B. Zeng, Y. He,
Free-matrix-based integral inequlaity for stability analysis of systems with
time-varying delay,  IEEE Transactions on Automatic Control, 60 (2015) 2768 - 2772.



\bibitem
{ZhH} X. Zhang, Y. Han, Y. Wang, C. Gong,
Abel lemma-based finite-sum inequality and its application to
stability analysis for linear discrete time-delay systems,
Automatica 57 (2015) 199-202.

\bibitem
{zha15} X. Zhang, Y. Han, Y. Wang, C. Gong,
Weighted orthogonal polynomial-based generalizaton of Wirtinger-type
integral inequalities for delayed continuous-time systems, (2015)
arXiv:1509.05085v1, math.OC 16 September.

\end{thebibliography}
\end{document}